\documentclass[12pt]{article}

\usepackage{amsmath,amsthm,amsfonts,amssymb,amscd,cite}
\usepackage{tikz}
\usepackage{color}
\usepackage{cite, url}

\newtheorem{theorem}{Theorem}[section]
\newtheorem{conjecture}[theorem]{Conjecture}

\begin{document}

\title{On the local metric dimension of $K_5$-free graphs}

\author{
Ali Ghalavand$^{a,}$\thanks{Corresponding author;  email:\texttt{alighalavand@nankai.edu.cn}}
\and 
Xueliang Li$^{a,b}$
}

\maketitle

\begin{center}
$^a$ Center for Combinatorics and LPMC, Nankai University, Tianjin 300071, China \\
\medskip
$^b$ School of Mathematical Sciences, Xinjiang Normal University, Urumchi, Xinjiang, China\\



\end{center}
\begin{center}
  Dedicated to Professor Fuji Zhang on the occasion of his 88th birthday.
\end{center}
\begin{abstract}
Let \( G \) be a graph with order \( n(G) \geq 5 \), local metric dimension \( \dim_l(G) \), and clique number \( \omega(G) \). In this paper, we investigate the local metric dimension of \( K_5 \)-free graphs and prove that \( \dim_l(G) \leq \lfloor\frac{2}{3}n(G)\rfloor \) when \( \omega(G) = 4 \). As a consequence of this finding, along with previous publications, we establish that if \( G \) is a \( K_5 \)-free graph, then \( \dim_l(G) \leq \lfloor\frac{2}{5}n(G)\rfloor \) when \( \omega(G) = 2 \), \( \dim_l(G) \leq \lfloor\frac{1}{2}n(G)\rfloor \) when \( \omega(G) = 3 \), and \( \dim_l(G) \leq \lfloor\frac{2}{3}n(G)\rfloor \) when \( \omega(G) = 4 \). Notably, these bounds are sharp for planar graphs. These results for graphs with a clique number less than or equal to 4 provide a positive answer to the conjecture stating that if \( n(G) \geq \omega(G) + 1 \geq 4 \), then \( \dim_l(G) \leq \left( \frac{\omega(G) - 2}{\omega(G) - 1} \right)n(G) \).
\end{abstract}
\noindent
\textbf{Keywords:} metric dimension; local metric dimension; clique number; planar graph

\medskip\noindent
\textbf{AMS Math.\ Subj.\ Class.\ (2020)}: 05C12, 05C69

\section{Introduction}
Let \( G \) be a simple connected graph with vertex set \( V(G) \) and edge set \( E(G) \). We denote the order of \( G \) as \( n(G) \) and its clique number as \( \omega(G) \). Let \( u \), \( v \), and \( w \) be arbitrary elements of \( V(G) \).
The notation \( N_G(u) \) represents the open neighborhood of \( u \), which is the set of vertices in \( V(G) \) that are adjacent to \( u \). The degree of vertex \( u \), denoted as \( d_G(u) \), refers to the number of vertices in \( N_G(u) \).
The notation \( d_G(u, v) \) denotes the distance between vertices \( u \) and \( v \), defined as the length of the shortest path in \( G \) connecting \( u \) and \( v \). We say that the vertices \( u \) and \( v \) are {\em distinguished} by \( w \), or equivalently, that \( w \) {\em distinguishes} \( u \) and \( v \), if \( d_G(u,w) \neq d_G(v,w) \).
For a positive integer \( k \), the notation \( [k] \) represents the set \( \{ 1, 2, \ldots, k \} \), and we define \( [0] = \emptyset \). For a positive integer \( n \), the notations \( K_n \) and \( P_n \) represent the complete graph and the path on \( n \) vertices, respectively. 
Let \( V' \) be a subset of \( V(G) \). The notation \( G[V'] \) denotes the induced subgraph of \( G \) with vertex set \( V' \), where two vertices are adjacent if and only if they are also adjacent in \( G \).
 Let \( G' \) be a distinct graph from \( G \). We denote the union of \( G \) and \( G' \) as \( G \cup G' \), which consists of the vertex set \( V(G) \cup V(G') \) and the edge set \( E(G) \cup E(G') \).
Now, consider \( H \) as an arbitrary subgraph of \( G \). The subgraph \( G - H \) is defined as the graph obtained by removing the vertices of \( H \) and the edges that have at least one endpoint in \( H \) from \( G \). If \( H' \) is another subgraph of \( G \), then \( E_G(H, H') \) denotes the collection of edges in \( G \) that connect one vertex in \( H \) to another vertex in \( H' \).

In this research, we focus on examining the local metric dimension of finite, simple, and connected graphs with a clique number of 4. First, let's define some key concepts. A {\em resolving set} for a graph \( G \) is a subset \( W \) of \( V(G) \) such that for any two distinct vertices \( u \) and \( v \) in \( V(G) - W \), there exists at least one vertex in \( W \) that can distinguish between \( u \) and \( v \). 
Similarly, a {\em local resolving set} of \( G \) is a subset \( W \) of \( V(G) \) such that for any adjacent vertices \( u \) and \( v \) in \( V(G) - W \), there is a vertex in \( W \) that can distinguish \( u \) from \( v \). The cardinalities of the smallest resolving sets and the smallest local resolving sets for \( G \) are referred to as the {\em metric dimension} \( \dim(G) \) and the {\em local metric dimension} \( \dim_l(G) \) of the graph \( G \), respectively. It is important to note that \( \dim_l(G) \leq \dim(G) \).

The concept of the metric dimension of graphs has an extensive history, originally defined by Harary and Melter ~\cite{13}, as well as Slater \cite{25}. Determining the metric dimension is known to be NP-complete for general graphs \cite{19}, and this complexity also extends to specific cases, such as planar graphs with a maximum degree of 6 \cite{6}. Research in this field is prolific, partly due to the metric dimension's wide range of real-world applications, which include robot navigation, image processing, privacy in social networks, and tracking intruders in networks. A 2023 overview \cite{survey2} of the essential results and applications of metric dimension includes over 200 references.

Research on the metric dimension has led to the exploration of various related concepts. A survey \cite{survey1} that focuses on these variations cites over 200 papers. One particularly interesting variant is the local metric dimension, introduced in 2010 by Okamoto et al. \cite{Okamoto1}. Similar to the standard metric dimension, the local metric dimension presents computational challenges \cite{9,10} and has been the subject of several studies \cite{Abrishami1, 3, 4,fitriani, Ghalavand1,klavzar-2023, lal-2023, 17,rodriguez-2016}, including research on the fractional local metric dimension \cite{javaid-2024}.

Okamoto et al. \cite{Okamoto1} established several significant relationships between the local metric dimension and the clique number: 
\begin{itemize}
\item  $\dim_l(G) = n(G) - 1$ if and only if $G \cong K_{n(G)}$;
\item $\dim_l(G) = n(G) - 2$ if and only if $\omega(G) = n(G) - 1$; 
\item $\dim_l(G) = 1$ if and only if $G$ is bipartite; 
\item $\dim_l(G) \geq \max\left \{ \lceil \log_2 \omega(G) \rceil, n(G) - 2^{n(G) - \omega(G)} \right\}$.
\end{itemize}
Additionally, Abrishami et al. \cite{Abrishami1} demonstrated that $\dim_l(G) \leq \frac{2}{5}n(G)$ when $\omega(G) = 2$ and $n(G) \geq 3$. Furthermore, one of the authors, along with others, proved in \cite{Ghalavand1} that \( \dim_l(G) \leq \left( \frac{\omega(G) - 1}{\omega(G)} \right) n(G) \), with equality occurring only if \( G \cong K_{n(G)} \). This result was first conjectured in \cite{Abrishami1}. The authors also presented the following conjecture.

\begin{conjecture} {\rm \cite[Conjecture 2]{Ghalavand1}}
\label{con2}
If $G$ is a graph with $n(G) \geq \omega(G)+1 \geq 4$, then 
\[\dim_l(G)\leq\left(\frac{\omega(G)-2}{\omega(G)-1}\right)n(G).\]
\end{conjecture}

It has been demonstrated in \cite{Ghalavand1} that if Conjecture~\ref{con2} is true, then the bound is asymptotically the best possible. 
Recently, the authors in \cite{Ghalavand2} confirmed this conjecture for all graphs with a clique number \(\omega(G)\) in the set \(\{n(G)-1, n(G)-2, n(G)-3\}\). They also characterized all graphs with order \(n\) and a local metric dimension of \(n-3\). Additionally, they established that:
\begin{itemize}
\item \(n(G)-4 \leq \dim_l(G) \leq n(G)-3\) when \(\omega(G) = n(G)-2\);
\item \(n(G)-8 \leq \dim_l(G) \leq n(G)-3\) when \(\omega(G) = n(G)-3\).
\end{itemize}
 Furthermore, the authors in \cite{Ghalavand3} confirmed this conjecture when \(\omega(G) = 3\).

In this paper, we present the following theorem, which positively addresses Conjecture \ref{con2} for graphs with a clique number of 4. It is important to note that Conjecture \ref{con2} remains unresolved for graphs \(G\) where \(5 \leq \omega(G) \leq n(G) - 4\).

\begin{theorem}\label{3th}
If $G$ is a  graph of order  $n(G)\geq5$ with clique number \(\omega(G) = 4\), then    $\dim_l(G) \leq\lfloor \frac{2}{3}n(G)\rfloor$. 
\end{theorem}

For any positive number \( t \), let \( G = tK_3 + K_1 \). This means that \( G \) is constructed from \( t \) disjoint complete graphs \( K_3 \) by adding a new vertex and connecting it to all vertices of the \( tK_3 \) components. It is easy to observe that \(\omega(G) = 4\) and $\dim_l(G) =\lfloor \frac{2}{3}n(G)\rfloor$. Therefore, there are infinitely many planar graphs \( G \) such that \(\omega(G) = 4\) and \(\dim_l(G) = \left\lfloor \frac{2}{3}n(G) \right\rfloor\).

\section{Proof of Theorem~\ref{3th}}
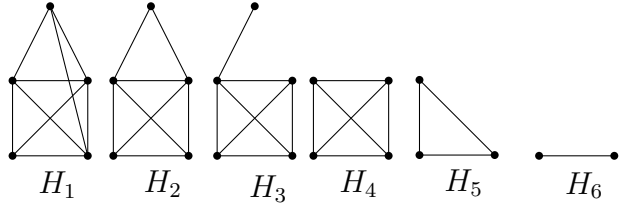
\begin{figure}[ht!]
\begin{center}
\begin{tikzpicture}
\clip(-7,1.2) rectangle (3,4.2);
\draw (-6,3)-- (-5,3);
\draw (-6,3)-- (-6,2);
\draw (-6,2)-- (-5,2);
\draw (-5,2)-- (-5,3);
\draw (-5,3)-- (-6,2);
\draw (-6,3)-- (-5,2);
\draw (-5.5,3.99)-- (-5,3);
\draw (-5.5,3.99)-- (-6,3);
\draw (-5.5,3.99)-- (-5,2);
\draw (-4.66,3)-- (-3.66,3);
\draw (-4.66,3)-- (-4.66,2);
\draw (-4.66,2)-- (-3.66,2);
\draw (-3.66,2)-- (-3.66,3);
\draw (-3.66,3)-- (-4.66,2);
\draw (-4.66,3)-- (-3.66,2);
\draw (-4.16,3.99)-- (-3.66,3);
\draw (-4.16,3.99)-- (-4.66,3);
\draw (-3.28,3)-- (-2.28,3);
\draw (-3.28,3)-- (-3.28,2);
\draw (-3.28,2)-- (-2.28,2);
\draw (-2.28,2)-- (-2.28,3);
\draw (-2.28,3)-- (-3.28,2);
\draw (-3.28,3)-- (-2.28,2);
\draw (-2.78,3.99)-- (-3.28,3);
\draw (-2,3)-- (-1,3);
\draw (-2,3)-- (-2,2);
\draw (-2,2)-- (-1,2);
\draw (-1,2)-- (-1,3);
\draw (-1,3)-- (-2,2);
\draw (-2,3)-- (-1,2);
\draw (-0.59,3.01)-- (-0.59,2.01);
\draw (-0.59,2.01)-- (0.41,2.01);
\draw (-0.59,3.01)-- (0.41,2.01);
\draw (1,2)-- (2,2);
\draw (-5.8,1.94) node[anchor=north west] {$H_1$};
\draw (-4.4,1.95) node[anchor=north west] {$H_2$};
\draw (-3.0,1.88) node[anchor=north west] {$H_3$};
\draw (-1.8,1.95) node[anchor=north west] {$H_4$};
\draw (-0.4,1.96) node[anchor=north west] {$H_5$};
\draw (1.2,1.93) node[anchor=north west] {$H_6$};
\begin{scriptsize}
\fill [color=black] (-6,3) circle (1.5pt);
\fill [color=black] (-5,3) circle (1.5pt);
\fill [color=black] (-6,2) circle (1.5pt);
\fill [color=black] (-5,2) circle (1.5pt);
\fill [color=black] (-5.5,3.99) circle (1.5pt);
\fill [color=black] (-4.66,3) circle (1.5pt);
\fill [color=black] (-3.66,3) circle (1.5pt);
\fill [color=black] (-4.66,2) circle (1.5pt);
\fill [color=black] (-3.66,2) circle (1.5pt);
\fill [color=black] (-4.16,3.99) circle (1.5pt);
\fill [color=black] (-3.28,3) circle (1.5pt);
\fill [color=black] (-2.28,3) circle (1.5pt);
\fill [color=black] (-3.28,2) circle (1.5pt);
\fill [color=black] (-2.28,2) circle (1.5pt);
\fill [color=black] (-2.78,3.99) circle (1.5pt);
\fill [color=black] (-2,3) circle (1.5pt);
\fill [color=black] (-1,3) circle (1.5pt);
\fill [color=black] (-2,2) circle (1.5pt);
\fill [color=black] (-1,2) circle (1.5pt);
\fill [color=black] (-0.59,3.01) circle (1.5pt);
\fill [color=black] (-0.59,2.01) circle (1.5pt);
\fill [color=black] (0.41,2.01) circle (1.5pt);
\fill [color=black] (1,2) circle (1.5pt);
\fill [color=black] (2,2) circle (1.5pt);
\end{scriptsize}
\end{tikzpicture}
\caption{ The graphs $H_1$, $H_2$, $\ldots$, $H_6$. }
\label{fig1}
\end{center}
\end{figure}

We begin the proof by outlining a crucial approach. First, let \( H_i \) (where \( i \in [6] \)) represent the graphs depicted in Fig.~\ref{fig1}. Let \( G \) be a graph with \( n(G) \geq 5 > \omega(G) \). In the following sections, we will sequentially identify maximum sets of vertex disjoint induced subgraphs within \( G \) and its induced subgraphs that are isomorphic to one of the \( H_i \) graphs. Although this selection is not necessarily unique, we will choose one specific selection and fix it for the purpose of this proof, ensuring that the following notation is well-defined.

\begin{itemize}
\item Let $\mathcal{H}_1(G)$ be a maximum set of vertex disjoint induced subgraphs of $G$ isomorphic to $H_1$.  
\item Set $G_1 = G$. For $i=2,3,\ldots,6$, let $\mathcal{H}_i(G)$ be a maximum set of vertex disjoint induced subgraphs of $G_{i}=G_{i-1} - \sum_{H \in \mathcal{H}_{i-1}(G)}H$ isomorphic to $H_i$.
\item Note that $G_{6}-\sum_{H \in \mathcal{H}_{6}(G)}H$ is a set of isolated vertices, let $\mathcal{H}_{7}(G)$ be the set of these induced subgraphs isomorphic to $K_1$.
\item For \( i \in [6] \), let  \(V_i= \cup_{H \in \mathcal{H}_{i}(G)} V(H) \).
\end{itemize}

It is clear that the sets \(V_i\) for \(i \in [7]\) form a partition of the vertex set of \(G\). Importantly, for \(i\) ranging from \(1\) to \(7\), the sets \(\mathcal{H}_i(G)\) may be empty. In the remainder of this text, we will consider the following conditions on the elements of \(\mathcal{H}_i(G)\).
\begin{itemize}
  \item For $i\in[6]$,  $\mathcal{H}_i(G)=\{H^i_1,\ldots,H^i_{|\mathcal{H}_i(G)|}\}$.
  \item For $i\in[|\mathcal{H}_1(G)|]$, $V(H^1_i)=\{h^1_{i_1},\ldots,h^1_{i_5}\}$, $d_{H^1_i}(h^1_{i_1})=d_{H^1_i}(h^1_{i_2})=d_{H^1_i}(h^1_{i_3})=4$, $d_{H^1_i}(h^1_{i_4})=d_{H^1_i}(h^1_{i_5})=3$.
   \item For $i\in[|\mathcal{H}_2(G)|]$, $V(H^2_i)=\{h^2_{i_1},\ldots,h^2_{i_5}\}$, $d_{H^2_i}(h^2_{i_1})=d_{H^2_i}(h^2_{i_2})=4$, $d_{H^2_i}(h^2_{i_3})=d_{H^2_i}(h^2_{i_4})=3$, and $d_{H^2_i}(h^2_{i_5})=2$.
    \item For $i\in[|\mathcal{H}_3(G)|]$, $V(H^3_i)=\{h^3_{i_1},\ldots,h^3_{i_5}\}$, $d_{H^3_i}(h^3_{i_1})=4$, $d_{H^3_i}(h^3_{i_2})=d_{H^3_i}(h^3_{i_3})=d_{H^3_i}(h^3_{i_4})=3$, and $d_{H^3_i}(h^3_{i_5})=1$.
    \item For $i\in[|\mathcal{H}_4(G)|]$, $V(H^4_i)=\{h^4_{i_1},\ldots,h^4_{i_4}\}$,  for $i\in[|\mathcal{H}_5(G)|]$, $V(H^5_i)=\{h^5_{i_1},\ldots,h^5_{i_3}\}$, and for $i\in[|\mathcal{H}_6(G)|]$, $V(H^6_i)=\{h^6_{i_1},h^6_{i_2}\}$.
\end{itemize}

Let $G_i$, where $i\in[6]$, be the subgraphs of $G$ defined above. By applying our assumptions that $\omega(G) \leq 4$ and the maximality of $\mathcal{H}_i(G)$ for $i \in [6]$, we can observe the following results.
\begin{enumerate}
  \item[(I)] If  $i\in[|\mathcal{H}_1(G)|]$, $j\in\{4,5\}$, and $v\in\,V(G-H^1_i)$, then $G[\{h^1_{i_1},h^1_{i_2},h^1_{i_3},h^1_{i_j},v\}]\not\cong\,K_5$. 
  \item[(II)] If  $i\in[|\mathcal{H}_2(G)|]$, $j\in\{2,3\}$, and $v\in\,V(G_2-H^2_i)$, then $G[\{h^2_{i_1},h^2_{i_4},h^2_{i_j},v\}]\not\cong\,K_4$. 
  \item[(III)] If  $i\in[|\mathcal{H}_3(G)|]$, $j\in\{1,2\}$, and $v\in\,V(G_3-H^3_i)$, then $G[\{h^3_{i_3},h^3_{i_4},h^3_{i_j},v\}]\not\cong\,K_4$. 
  \item[(IV)] If  $i\in[|\mathcal{H}_4(G)|]$, $j\in[4]$, and $v\in\,V(G_4-H^4_i)$, then $h^4_{i_j}v\not\in\,E(G)$. 
  \item[(V)] If $i\in[|\mathcal{H}_4(G)|]$, $j\in[4]$, $v\in\,V(G-G_4)$, and $h^4_{i_j}v\in\,E(G)$, then  there exists an element  $l\in([4]-\{j\})$ such that $G[h^4_{i_{j}},h^4_{i_{l}},v]\cong\,P_3$.
  \item[(VI)] If  $i\in[|\mathcal{H}_5(G)|]$ and $v\in\,V(G_5-H^5_i)$, then $G[V(H^5_i)\cup\{v\}]\not\cong\,K_4$.
  \item[(VII)] If  $i\in[|\mathcal{H}_6(G)|]$ and $v\in\,V(G_6-H^6_i)$, then $G[V(H^6_i)\cup\{v\}]\not\cong\,K_3$.
\end{enumerate}

In the following processes, we will construct a set \( S \) such that \( |S| \leq \frac{2}{3} n(G) \), ensuring that \( S \) remains a local resolving set for \( G \). We start with \( S = \emptyset \). To move forward, we need to introduce an additional notation. Let \( \mathcal{X} \subseteq \mathcal{H}_4(G) \) and  \( H \in \mathcal{H}_i(G) \), where \( i\in [3] \). The notation \( \tau_{i,4}(H,\mathcal{X}) \) represents the set of elements \( X \) in \( \mathcal{X} \) such that \( E_G(H,X) \neq \emptyset \).

\medskip\noindent
\underline{$1^{\rm st}$ process:}
\begin{enumerate}
\item[(1.1)] Set $S = \emptyset$, $i=1$, and $\mathcal{X}=\mathcal{H}_4(G)$, and then go to (1.2).
\item[(1.2)] If $i>|\mathcal{H}_1(G)|$, then return  $S$ and $\mathcal{X}$,  and end the process,  otherwise go to (1.3).  
\item[(1.3)] If  $|\tau_{1,4}(H_i^1,\mathcal{X})|=0$, then set
\begin{align*}
S & = S\cup\,(V(H_i^1)-\{h^1_{i_4},h^1_{i_5}\}),\\
i & =i+1,
\end{align*}
and proceed to (1.2), otherwise go to (1.4).
\item[(1.4)] If  $|\tau_{1,4}(H_i^1,\mathcal{X})|=1$ and $\tau_{1,4}(H_i^1,\mathcal{X})=\{X_1\}$, then choose two distinct elements $l_1,l_2\in[5]$ and two distinct elements $x^1_1,x^1_2\in\,V(X_1)$ such that $G[\{h^1_{i_{l_1}},x^1_1,x^1_2\}]\cong\,P_3$, and then set
\begin{align*}
S & = S\cup\,(V(H_i^1)-\{h^1_{i_{l_2}}\})\cup(V(X_1)-\{x^1_1,x^1_2\}),\\
i & =i+1,\\
\mathcal{X}&=\mathcal{X}-\tau_{1,4}(H_i^1,\mathcal{X}),
\end{align*}
and proceed to (1.2), otherwise go to (1.5).
\item[(1.5)] If $|\tau_{1,4}(H_i^1,\mathcal{X})|=2$ and $\tau_{1,4}(H_i^1,\mathcal{X})=\{X_1,X_2\}$, then for an element $l\in[5]$, for which $G[(V(H_i^1)-\{h^1_{i_l}\})\cup V(X_1)\cup V(X_2)]$ is connected, choose elements $x^1_1,x^1_2\in V(X_1)$ and $x^2_1,x^2_2\in V(X_2)$ such that for $z_1,z_2\in([5]-\{l\})$, not necessarily distinct, we have $G[\{h^1_{i_{z_1}},x^1_1,x^1_2\}]\cong\,P_3$ and $G[\{h^1_{i_{z_2}},x^2_1,x^2_2\}]\cong\,P_3$, and then set
\begin{align*}
S & = S\cup(V(H_i^1)-\{h^1_{i_l}\})\cup(V(X_1)-\{x^1_1,x^1_2\})\cup(V(X_2)-\{x^2_1,x^2_2\}),\\
i & =i+1,\\
\mathcal{X}&=\mathcal{X}-\tau_{1,4}(H_i^1,\mathcal{X}),
\end{align*}
and proceed to (1.2), otherwise go to (1.6).
\item[(1.6)] If  $|\tau_{1,4}(H_i^1,\mathcal{X})|\geq3$ and $\tau_{1,4}(H_i^1,\mathcal{X})=\{X_1,\ldots,X_{|\tau_{1,4}(H_i^1,\mathcal{X})|}\}$, then for $l$ from $1$ to $|\tau_{1,4}(H_i^1,\mathcal{X})|$, choose two distinct elements $x^l_1,x^l_2 \in V(X_l)$ such that for an element $z \in [5]$, $G[\{h^1_{i_z},x^l_1,x^l_2\}]\cong\,P_3$, and then set
\begin{align*}
S & = S\cup\,V(H_i^1)\cup_{l=1}^{|\tau_{1,4}(H_i^1,\mathcal{X})|}(V(X_l)-\{x^l_1,x^l_2\}),\\
i & =i+1,\\
\mathcal{X}&=\mathcal{X}-\tau_{1,4}(H_i^1,\mathcal{X}),
\end{align*}
and proceed to (1.2).
\end{enumerate}

\noindent
\underline{$2^{\rm nd}$ process:} 
\begin{enumerate} 
\item[(2.1)] Consider the sets $S$ and $\mathcal{X}$ that are returned in the $1^{\rm st}$ process, and then set $i=1$, and go to (2.2).
\item[(2.2)] If $i>|\mathcal{H}_2(G)|$, then return  $S$ and $\mathcal{X}$,  and end the process,  otherwise go to (2.3).  
\item[(2.3)] If  $|\tau_{2,4}(H_i^2,\mathcal{X})|=0$, then set
\begin{align*}
S & = S\cup\,(V(H_i^2)-\{h^2_{i_2},h^2_{i_3}\}),\\
i & =i+1,
\end{align*}
and proceed to (2.2), otherwise go to (2.4).
\item[(2.4)] If \( |\tau_{2,4}(H_i^2, \mathcal{X})| = 1 \) and \( \tau_{2,4}(H_i^2, \mathcal{X}) = \{X_1\} \), then choose an element \( l_1 \in [4] \) such that \( G[(V(H_i^2) - \{h^2_{i_{l_1}}\}) \cup V(X_1)] \) is connected. Additionally, select two distinct elements \( x^1_1, x^1_2 \in V(X_1) \) such that for an element \( l_2 \in ([5] - \{l_1\}) \), it holds that $G[\{h^2_{i_{l_2}},x^1_1,x^1_2\}]\cong\,P_3$. Then set
\begin{align*}
S & = S\cup\,(V(H_i^2)-\{h^2_{i_{l_1}}\})\cup(V(X_1)-\{x^1_1,x^1_2\}),\\
i & =i+1,\\
\mathcal{X}&=\mathcal{X}-\tau_{2,4}(H_i^2,\mathcal{X}),
\end{align*}
and proceed to (2.2), otherwise go to (2.5).
\item[(2.5)] If $|\tau_{2,4}(H_i^2,\mathcal{X})|=2$ and $\tau_{2,4}(H_i^2,\mathcal{X})=\{X_1,X_2\}$, then for an element $l_1\in[4]$, for which $G[(V(H_i^2)-\{h^2_{i_{l_1}}\})\cup V(X_1)\cup V(X_2)]$ is connected, choose elements $x^1_1,x^1_2\in V(X_1)$ and $x^2_1,x^2_2\in V(X_2)$ such that for $l_2,l_3\in([5]-\{l_1\})$, not necessarily distinct, we have $G[\{h^2_{i_{l_2}},x^1_1,x^1_2\}]\cong\,P_3$ and $G[\{h^2_{i_{l_3}},x^2_1,x^2_2\}]\cong\,P_3$, and then set
\begin{align*}
S & = S\cup(V(H_i^2)-\{h^2_{i_{l_1}}\})\cup(V(X_1)-\{x^1_1,x^1_2\})\cup(V(X_2)-\{x^2_1,x^2_2\}),\\
i & =i+1,\\
\mathcal{X}&=\mathcal{X}-\tau_{2,4}(H_i^2,\mathcal{X}),
\end{align*}
and proceed to (2.2), otherwise go to (2.6).
\item[(2.6)] If  $|\tau_{2,4}(H_i^2,\mathcal{X})|\geq3$ and $\tau_{2,4}(H_i^2,\mathcal{X})=\{X_1,\ldots,X_{|\tau_{2,4}(H_i^2,\mathcal{X})|}\}$, then for $l$ from $1$ to $|\tau_{2,4}(H_i^2,\mathcal{X})|$, choose two distinct elements $x^l_1,x^l_2 \in V(X_l)$ such that for an element $z \in [5]$, $G[\{h^2_{i_z},x^l_1,x^l_2\}] \cong\,P_3$, and then set
\begin{align*}
S & = S\cup\,V(H_i^2)\cup_{l=1}^{|\tau_{2,4}(H_i^2,\mathcal{X})|}(V(X_l)-\{x^l_1,x^l_2\}),\\
i & =i+1,\\
\mathcal{X}&=\mathcal{X}-\tau_{2,4}(H_i^2,\mathcal{X}),
\end{align*}
and proceed to (2.2).
\end{enumerate}

\noindent
\underline{$3^{\rm rd}$ process:} 
\begin{enumerate} 
\item[(3.1)] Consider the sets $S$ and $\mathcal{X}$ that are returned in the $2^{\rm nd}$ process, and then set $i=1$, and go to (3.2).
\item[(3.2)] If $i>|\mathcal{H}_3(G)|$, then return  $S$,  and end the process,  otherwise go to (3.3).  
\item[(3.3)] If  $|\tau_{3,4}(H_i^3,\mathcal{X})|=0$, then set
\begin{align*}
S & = S\cup\,(V(H_i^3)-\{h^3_{i_1},h^3_{i_2}\}),\\
i & =i+1,
\end{align*}
and proceed to (3.2), otherwise go to (3.4).
\item[(3.4)] If \( |\tau_{3,4}(H_i^3, \mathcal{X})| = 1 \) and \( \tau_{3,4}(H_i^2, \mathcal{X}) = \{X_1\} \), then choose an element \( l_1 \in [4] \) such that \(E_G( G[V(H_i^3) - \{h^3_{i_{l_1}}\}],X_1)\neq\emptyset \). Additionally, select two distinct elements \( x^1_1, x^1_2 \in V(X_1) \) such that for an element \( l_2 \in ([5] - \{l_1\}) \), it holds that
     $G[\{h^3_{i_{l_2}},x^1_1,x^1_2\}]\cong\,P_3$. Then set
\begin{align*}
S & = S\cup\,(V(H_i^3)-\{h^3_{i_{l_1}}\})\cup(V(X_1)-\{x^1_1,x^1_2\}),\\
i & =i+1,\\
\mathcal{X}&=\mathcal{X}-\tau_{3,4}(H_i^3,\mathcal{X}),
\end{align*}
and proceed to (3.2), otherwise go to (3.5).
\item[(3.5)] If $|\tau_{3,4}(H_i^3,\mathcal{X})|=2$ and $\tau_{3,4}(H_i^3,\mathcal{X})=\{X_1,X_2\}$, then for an element $l_1\in[4]$, for which $E_G(G[V(H_i^3)-\{h^3_{i_{l_1}}\}],X_1)\neq\emptyset$ and $E_G(G[V(H_i^3)-\{h^3_{i_{l_1}}\}],X_2)\neq\emptyset$, choose elements $x^1_1,x^1_2\in V(X_1)$ and $x^2_1,x^2_2\in V(X_2)$ such that for $l_2,l_3\in([5]-\{l_1\})$, not necessarily distinct, we have $G[\{h^3_{i_{l_2}},x^1_1,x^1_2\}]\cong\,P_3$ and $G[\{h^3_{i_{l_3}},x^2_1,x^2_2\}]\cong\,P_3$, and then set
\begin{align*}
S & = S\cup(V(H_i^3)-\{h^3_{i_{l_1}}\})\cup(V(X_1)-\{x^1_1,x^1_2\})\cup(V(X_2)-\{x^2_1,x^2_2\}),\\
i & =i+1,\\
\mathcal{X}&=\mathcal{X}-\tau_{3,4}(H_i^3,\mathcal{X}),
\end{align*}
and proceed to (3.2), otherwise go to (3.6).
\item[(3.6)] If  $|\tau_{3,4}(H_i^3,\mathcal{X})|\geq3$ and $\tau_{3,4}(H_i^3,\mathcal{X})=\{X_1,\ldots,X_{|\tau_{3,4}(H_i^3,\mathcal{X})|}\}$, then for $l$ from $1$ to $|\tau_{3,4}(H_i^3,\mathcal{X})|$, choose two distinct elements $x^l_1,x^l_2 \in V(X_l)$ such that for an element $z \in [5]$, $G[\{h^3_{i_z},x^l_1,x^l_2\}]\cong\,P_3$, and then set
\begin{align*}
S & = S\cup\,V(H_i^3)\cup_{l=1}^{|\tau_{3,4}(H_i^3,\mathcal{X})|}(V(X_l)-\{x^l_1,x^l_2\}),\\
i & =i+1,\\
\mathcal{X}&=\mathcal{X}-\tau_{3,4}(H_i^3,\mathcal{X}),
\end{align*}
and proceed to (3.2).
\end{enumerate}

\noindent
\underline{$4^{\rm th}$ process:} 
\begin{enumerate} 
\item[(4.1)]  Consider the set \( S \) that is obtained in the $3^{\rm rd}$ process, and then set \( i = 1 \) before proceeding to step (4.2).
\item[(4.2)] If $i>|\mathcal{H}_5(G)|$, then return  $S$,  and end the process,  otherwise go to (4.3).  
\item[(4.3)] Set 
\begin{align*}
S & = S\cup\{h^5_{i_1},h^5_{i_2}\},\\
i & =i+1,
\end{align*}
and proceed to (4.2).
\end{enumerate}

\noindent
\underline{$5^{\rm th}$ process:} 
\begin{enumerate} 
\item[(5.1)] Consider the set \( S \) that is obtained in the $4^{\rm th}$ process, and then set \( i = 1 \) before proceeding to step (5.2).
\item[(5.2)] If \( i > |\mathcal{H}_6(G)| \), then set \( S = S \cup \mathcal{H}_7(G) \), return \( S \), and end the process. Otherwise, proceed to (5.3). 
\item[(5.3)] Set 
\begin{align*}
S & = S\cup\{h^6_{i_1}\},\\
i & =i+1,
\end{align*}
and proceed to (5.2).
\end{enumerate}

Now, let \( S \) be the subset of \( V(G) \) obtained from the five processes described above. For an integer \( l \), define \( \lambda = \frac{2}{3}(4l + 5) \). It is evident that \( \lambda \geq 3 \) when \( l = 0 \), \( \lambda = 6 \) when \( l = 1 \), \( \lambda > 8 \) when \( l = 2 \), and \( \lambda \geq 2l + 5 \) when \( l \geq 3 \). 
Additionally, we have that \( \frac{2}{3} \cdot 4 > 2 \), \( \frac{2}{3} \cdot 3 = 2 \), and \( \frac{2}{3} \cdot 2 > 1 \). Thus, by applying the methods used to create \( S \), we can observe that \( |S| \leq \frac{2}{3}n(G) \).
Moreover, by utilizing the cases (I) to (VII) mentioned earlier, it follows that \( S \) serves as a local resolving set for \( G \). Since \( \dim_l(G) \) is an integer, we have proved Theorem~\ref{3th}.

\section*{Acknowledgments}

We would like to thank the anonymous referees for their valuable comments and suggestions. The research  was supported by the NSFC No.\ 12131013. 

\section*{Conflicts of interest} 

The authors declare no conflict of interest.

\section*{Data availability} 

No data was used in this investigation.


\begin{thebibliography}{99}

\bibitem{Abrishami1} 
G. Abrishami, M. A. Henning, M. Tavakoli,  
Local metric dimension for graphs with small clique numbers,
Discrete Math.\ 345 (2022) Paper 112763.

\bibitem{3} 
G.A. Barrag\'{a}n-Ram\'{\i}rez, A. Estrada-Moreno, Y. Ram\'{\i}rez-Cruz, J.A. Rodr\'{\i}guez-Vel\'{a}zquez, The local metric dimension of the lexicographic product of graphs, 
Bull.\ Malays.\ Math.\ Sci.\ Soc.\ 42 (2019) 2481--2496.

\bibitem{4} 
G.A. Barrag\'{a}n-Ram\'{\i}rez, J.A. Rodr\'{\i}guez-Vel\'{a}zquez, 
The local metric dimension of strong product graphs, 
Graphs Combin.\ 32 (2016) 1263--1278.

\bibitem{6} 
J. Diaz, O. Pottonen, M. Serna, E.J. van Leeuwen, 
Complexity of metric dimension on planar graphs, 
J.\ Comput.\ Syst.\ Sci.\ 83 (2017) 132--158.

\bibitem{9} 
H. Fernau, J.A. Rodr\'{\i}guez-Vel\'{a}zquez, 
On the (adjacency) metric dimension of corona and strong product graphs and their local variants, combinatorial and computational results, 
Discrete Appl.\ Math. 236 (2018) 183--202.

\bibitem{10} 
H. Fernau, J.A. Rodr\'{\i}guez-Vel\'{a}zquez, 
Notions of metric dimension of corona products: combinatorial and computational results, 
Lecture Notes Comput.\ Sci.\ 8476 (2014) 153--166.

\bibitem{fitriani} 
D.~Fitriani, S.W.~Saputro, 
The local metric dimension of amalgamation of graphs,
Electron.\ J.\ Graph Theory Appl.\ (EJGTA) 12 (2024) 125--146.

\bibitem{Ghalavand1} 
A. Ghalavand, M. A. Henning, M. Tavakoli, 
On a conjecture about the local metric dimension of graphs,
Graphs Combin.\ 39 (2023) Paper 5.

\bibitem{Ghalavand2} 
A. Ghalavand, S. Klav\v zar, X. Li, 
Interplay between the local metric dimension and the clique number of a graph, 
\url{arXiv:2412.17074} [math.CO] (22 Dec 2024). 

\bibitem{Ghalavand3} 
A. Ghalavand, S. Klav\v zar, X. Li, 
On the local metric dimension of $K_4$-free graphs, 
\url{arXiv:2506.00414} [math.CO] (31 May 2025). 

\bibitem{13} 
F. Harary, R.A. Melter, 
The metric dimension of a graph, 
Ars Combin.\ 2 (1976) 191--195.

\bibitem{javaid-2024}
I. Javaid, H. Benish, M. Murtaza, 
The fractional local metric dimension of graphs,
Contrib.\ Discrete Math.\ 19 (2024) 163--177.

\bibitem{19} 
S. Khuller, B. Raghavachari, A. Rosenfeld, 
Landmarks in graphs, 
Discrete Appl.\ Math.\ 70 (1996) 217--229.

\bibitem{klavzar-2023} 
S. Klav\v{z}ar, D. Kuziak,
Nonlocal metric dimension of graphs,
Bull.\ Malays.\ Math.\ Sci.\ Soc.\ 46 (2023) Paper 66.

\bibitem{17} 
S. Klav\v{z}ar, M. Tavakoli, 
Local metric dimension of graphs: generalized hierarchical products and some applications, 
Appl.\ Math.\ Comput.\ 364 (2020) Paper 124676.

\bibitem{survey1}
D.~Kuziak, I.G.~Yero,
Metric dimension related parameters in graphs: A survey on combinatorial, computational and applied results,  
\url{arXiv:2107.04877} [math.CO] (10 Jul 2021). 

\bibitem{lal-2023}
S. Lal, V.K. Bhat, 
On the local metric dimension of generalized wheel graph,
Asian-Eur.\ J.\ Math.\ 16 (2023) Paper 2350194.

\bibitem{Okamoto1} 
F. Okamoto, B. Phinezy, P. Zhang,  
The local metric dimension of a graph, 
Math.\ Bohem.\ 135 (2010) 239--255.

\bibitem{rodriguez-2016}
J.A.~Rodr\'iguez-Vel\'azquez, G.A.~Barrag\'an-Ram\'irez, C.~Garc\'ia G\'omez, 
On the local metric dimension of {C}orona product graphs,
Bull.\ Malays.\ Math.\ Sci.\ Soc.\ 39 (2016) S157--S173.

\bibitem{25} 
P.J. Slater, 
Leaves of trees, 
Congress.\ Numer.\ 14 (1975) 549--559. 

\bibitem{survey2}
R.~C.~Tillquist, R.~M.~Frongillo, M.~E.~Lladser. 
Getting the lay of the land in discrete space: a survey of metric dimension and its applications. 
SIAM Rev.\ 65 (2023) 919--962.

\end{thebibliography}
\end{document}